\documentclass{article}
\usepackage{fullpage}
\usepackage{graphics,psfrag,epsfig}
\usepackage{amsfonts,latexsym,eucal,amsmath,amsthm,amssymb} 

\newtheorem{assumption}{Assumption}
\newtheorem{theorem}{Theorem}
\newtheorem{lemma}{Lemma}

\newcommand{\RR}{\mathbb{R}}

\newcommand{\ZZ}{\mathbb{Z}}

\newcommand{\ES}{\mathcal{E}}
\newcommand{\Ltwo}{{L^2}}
\newcommand{\LL}[1]{\lvert#1\rvert_{{\Ltwo}}}

\newcommand{\hf}{\frac{1}{2}}

\newcommand{\kb}{\bar{k}}

\newcommand{\problem}{$\mbox{\ref{eq:NSFourierVorticity}}^\alpha_{\mathcal{Z}}$}
\newcommand{\problemB}{$\mbox{\ref{3DvorticityFourier2}}^\alpha_{\mathcal{Z}}$}
\newcommand{\problemBr}{${\mbox{\ref{3DvorticityFourier2}}^\alpha_{\mathcal{Z}}}^{(1)}$}
\newcommand{\problemBi}{${\mbox{\ref{3DvorticityFourier2}}^\alpha_{\mathcal{Z}}}^{(2)}$}
\title{An Elementary Proof of the Existence and Uniqueness Theorem for
  the Navier-Stokes Equations}
\author{ J. C. Mattingly\footnote{Department of Mathematics, Stanford
    University, Stanford CA 94305.} and Ya.\ G. Sinai\footnote{Department of Mathematics,
    Princeton University, Princeton NJ 08544.} }
\newcommand{\pd}{\partial}
\date{ November 19, 1998, Revised March 10th, 1999}
\begin{document}

\maketitle

\section{ Introduction}
The purpose of this paper is to show that some results concerning
solutions of the Navier-Stokes systems can be proven by purely
elementary methods. In two-dimensions with periodic boundary conditions,
the Navier-Stokes system has the form
\begin{align}  
\label{eq:NS}
  \frac{\partial u_1}{\pd t} + u_1 \frac{\partial u_1}{\pd x_1} + u_2
  \frac{\partial u_1}{\pd x_2} &= \nu \Delta u_1 - \frac{\partial
  p}{\partial x_1} + f_1(x_1,x_2,t) \ , \\
 \frac{\partial u_2}{\pd t} + u_1 \frac{\partial u_2}{\pd x_1} + u_2
  \frac{\partial u_2}{\pd x_2} &= \nu \Delta u_2 - \frac{\partial
  p}{\partial x_2} + f_2(x_1,x_2,t)  \ , \notag \\ 
\frac{\partial u_1}{\partial x_1}+\frac{\partial u_2}{\partial x_2}&=0
  \ . \notag
\end{align}

Here $\nu$ is the viscosity, $p$ is the pressure, and $f_1$, $f_2$ are the
components of an external forcing which may be time-dependent. As our
setting is periodic, the functions $u_1$, $u_2$, $\nabla p$, $f_1$, and
$f_2$ are all periodic in $x$. For simplicity, we take the period to
be one.

The first existence and uniqueness theorems for weak solutions of
\eqref{eq:NS} were proven by Leray (\cite{b:Leray34}) in whole plane
$\RR^2$. Later these results were extended by E. Hopf (see
\cite{b:Hopf51}). In 1962, Ladyzenskaya proved existence and uniqueness
results for strong solutions for general two-dimensional domains
\cite{b:Ladyzhenskaya69}. V. Yudovich, C. Foias, R. Teman, P.
Constantin, and others developed strong methods which provided deep
insights into the dynamics described by \eqref{eq:NS} (see
\cite{b:Yudovich89,b:Temam79,b:Temam95,b:CoFo88}).

The purpose of this paper is to present elementary proofs of three
theorems. These theorems imply the existence and uniqueness of smooth solutions
of \eqref{eq:NS} and shed some additional light on the dissipative character
of the dynamics. We will also discuss what our techniques can give in
the three-dimensional setting.

In two-dimensions, it is useful to consider the vorticity
$\omega(x_1,x_2,t) = \frac{\partial u_1(x_1,x_2,t)}{\partial x_2} -
\frac{\partial u_2(x_1,x_2,t)}{\partial x_1}$. The equation governing
$\omega$ has the form ( see \cite{b:ChorinMarsden93,b:DoGi95} )
\begin{align}
\label{eq:NSvorticity}
  \frac{\partial \omega}{\partial t} + u_1 \frac{\partial
    \omega}{\partial x_1}+ u_2 \frac{\partial \omega}{\partial x_2} =
    \nu \Delta \omega + g(x_1,x_2,t)
\end{align}
where $g(x_1,x_2,t) = \frac{\partial f_1(x_1,x_2,t)}{\partial x_2} -
\frac{\partial f_2(x_1,x_2,t)}{\partial x_1}$.  We will need
$g(x_1,x_2)$ to posses a modicum of spatial smoothness; this will be
made precise shortly.

In our two-dimensional setting, the systems \eqref{eq:NS} and
\eqref{eq:NSvorticity} are equivalent.
Expanding $\omega$ in Fourier series where $\omega(x_1,x_2,t) = \sum_{k \in
  \ZZ^2 } \omega_k(t) e^{2 \pi i (x,k)}$ with $x=(x_1,x_2)$, we can
write a coupled ODE-system for the modes $\omega_k(t)$ (see
\cite{b:DoGi95}).
\begin{align}
  \label{eq:NSFourierVorticity}
  \frac{d \omega_k}{d t} + 2 \pi i \sum_{l_1+l_2=k} \omega_{l_1}
  \omega_{l_2} \frac{(k,l_2^\perp)}{(l_2,l_2)}= -4 \pi^2 \nu \lvert k
  \rvert^2 \omega_k + g_k(t)
\end{align}
where $k \in \ZZ^2$, $|k|=\sqrt{k_1^2+k_2^2}$,
$l^\perp=(l^{(1)},l^{(2)})^\perp=(-l^{(2)},l^{(1)})$, and $g_k(t)$ are the
spatial Fourier modes of the function $g(x,t)$. Since $\omega$ is
real, we know $\omega_{-k}=\bar{\omega}_k$.  Furthermore, we
always assume that $\omega_0=0$. The system
\eqref{eq:NSFourierVorticity} is the Galerkin system corresponding to
\eqref{eq:NSvorticity}. A finite dimensional approximation of this
Galerkin system can be associated to any finite subset $\mathcal{Z}$ of
$\ZZ^2$ by setting $\omega_k(t)=0$ for all $k$ outside of
$\mathcal{Z}$. In the following, we will implicitly assume that
$\mathcal{Z}$ is centrally-symmetric, that is  if $k\in  \mathcal{Z}$
then  $-k\in  \mathcal{Z}$.

 In fact, we will study a slightly more general version of
\eqref{eq:NSvorticity} where the Laplacian is replaced by the operator
$|\nabla|^\alpha$ with $\alpha >1$. This leads to a version of
\eqref{eq:NSFourierVorticity} which we index by the choice of $\alpha$
and by the finite index set $\mathcal{Z}$, $\mathcal{Z} \subset
\ZZ^2$, indicating which modes are included in the Galerkin approximation.
In short, we consider the finite dimensional ODE system
\begin{align*}
  \frac{d \omega_k}{d t} + 2 \pi i \sum_{\substack{l_1+l_2=k\\ l_1,l_2
  \in \mathcal{Z}}} \omega_{l_1}
  \omega_{l_2} \frac{(k,l_2^\perp)}{(l_2,l_2)}= -4 \pi^2 \nu \lvert k
  \rvert^\alpha \omega_k + g_k  \ .\tag{\problem}
\end{align*}

We now state the assumptions on the coefficients $g_k(t)$ to be used
at various times during our discussion.
\begin{assumption} 
\label{a:L2Finite}
The forcing  $f(x,t)=(f_1(x,t),f_2(x,t))$ is such
  that ${\nobreak g^*=\sup_{t \in[0,\infty)} \LL{g(\cdot,t)} < \infty}$.
\end{assumption}

\begin{assumption}
\label{a:powerLaw}
  For some $r$, there exists a constant $\mathcal{G}(r)>0$
  such that 
  \begin{align*}
    \sup_{t \in [0,\infty)} |g_k(t)| \leq \frac{\mathcal{G}(r)}{|k|^{r-\alpha+\epsilon}}
  \end{align*}
for some $\epsilon >0$ and all $k \in \ZZ^2\backslash 0$. The constant $\alpha$ is the same
as in (\problem).
\end{assumption}

\begin{assumption}
\label{a:exponential}
For some $r$ and $\gamma > 0$, there exists a constant
$\mathcal{G}(r,\gamma)>0$ such that
 \begin{align*}
   \sup_{t \in [0,\infty)} |g_k(t)| \leq
   \frac{\mathcal{G}(r,\gamma)}{|k|^{r-\alpha+\epsilon}}e^{-\gamma|k|^{1+\delta}}
  \end{align*}
  for some $\delta >0$, $\epsilon >0$, and all $k \in \ZZ^2\backslash 0$. Again,
  the constant $\alpha$ is the same as in (\problem).
\end{assumption}
Observe that assumption \ref{a:exponential} implies assumption
\ref{a:powerLaw}.  Critical to our discussion is that for (\problem)
we have the so-called enstrophy estimate. Namely, if $ \ES(0)=\int
\omega^2(x_1,x_2,0) dx_1 dx_2 = \sum_k \lvert \omega_k(0) \rvert^2 <
\infty$ then one can find $\ES^*$ depending only on $ \ES(0)$, $\nu$,
$\sup_{t \in [0,\infty)}\LL{g(\cdot,t)}$, and $\alpha$ such that
${\ES(t)=\int \omega^2(x_1,x_2,t) dx_1 dx_2 \leq \ES^*}$ for all solutions to
(\problem).  It is important to
note that $\ES^*$ is independent of the set $\mathcal{Z}$ which defines
the Galerkin approximation. This enstrophy estimate holds if the
forcing satisfies assumption \ref{a:L2Finite} (see e.g.
\cite{b:CoFo88, b:DoGi95,b:Temam79}).

Now we are ready to formulate our theorems.

\begin{theorem}
  \label{thm:algebraicAlgebraic}
  Assume the forcing satisfies  assumption \ref{a:L2Finite} and
  \ref{a:powerLaw} for some $r>1$ and $\mathcal{G}(r)>0$. If for some $\mathcal{D}_1 < \infty$ 
  \begin{align*}
    \lvert \omega_k(0) \rvert \leq \frac{\mathcal{D}_1}{|k|^r}
  \end{align*}
  then one can find a $\mathcal{D}'_1 <\infty$, depending only on
  $\mathcal{D}_1$, $\nu$, $g^*$, and $\mathcal{G}$, such that any
  solution to (\problem) with these initial conditions satisfies
  \begin{align*}
    |\omega_k(t)| \leq \frac{\mathcal{D}_1'}{|k|^r}
  \end{align*}
  for all $t >0$. In particular, $\mathcal{D}_1'$ is independent of
  the set $\mathcal{Z}$ defining the Galerkin approximation.
\end{theorem}

An existence and uniqueness theorem for \eqref{eq:NSFourierVorticity}
follows from theorem \ref{thm:algebraicAlgebraic} by now standard
considerations (see \cite{b:CoFo88,b:DoGi95,b:Temam79}). We briefly
recall the general line of the argument. By the Sobolev embedding
theorem, the Galerkin approximations are trapped in a compact subset
of $L^2$ of the 2-torus. This guarantees the existence of a limit
point which can be shown to satisfy \eqref{eq:NSFourierVorticity}.
Using the the regularity inherited from the Galerkin approximations,
one then shows that there is a unique solution to
\eqref{eq:NSFourierVorticity}. Gallavotti \cite{b:Gallavotti96}
contains a similar proof of a similar statement.


\begin{theorem}
  \label{thm:exponentialExponential}
  Assume that assumption \ref{a:exponential} holds for some $r>1$,
  $\gamma >0$, and $\mathcal{G}(r,\gamma)>0$. If the initial
  conditions satisfy
  \begin{align*}
    |\omega_k(0)| \leq \frac{\mathcal{D}_2}{|k|^r}e^{-\gamma_2 |k|}
  \end{align*}
  for some $\mathcal{D}_2 < \infty$ and $\gamma_2>0$, then one can find
  a $\mathcal{D}'_2< \infty$ and a $\gamma_2'>0$, depending only on $\mathcal{D}_2$,
  $\gamma$, $r$, $\nu$, $g^*,\mathcal{G}$, such that any solution to (\problem)
  starting from these initial conditions satisfies
    \begin{align*}
      |\omega_k(t)| \leq \frac{\mathcal{D}'_2}{|k|^r} e^{-\gamma'_2 |k|}
    \end{align*} 
    for all $t > 0$. In particular, the constants $\mathcal{D}'_2$ and
    $\gamma'$ are independent of the set $\mathcal{Z}$ defining the
    Galerkin approximation.
  \end{theorem}

Theorem \ref{thm:exponentialExponential} shows that equation
\eqref{eq:NSvorticity} preserves the class of real analytic functions
on the 2-torus.
\begin{theorem}
  \label{thm:algebraicExponential}
  Assume that assumption \ref{a:exponential} holds for some $r>1$,
  $\gamma >0$, and $\mathcal{G}(r,\gamma)>0$.  If the initial conditions satisfy
  \begin{align*}
    |\omega_k(0)| \leq \frac{\mathcal{D}_3}{|k|^r}
  \end{align*}
  then for any $t_0>0$, one can find a $\mathcal{D}_3' >0$ and
  a $\gamma_3' > 0$ such that any solution to (\problem) with these
  initial conditions satisfies
  \begin{align*}
    |\omega_k(t_0)| \leq \frac{\mathcal{D}_3'}{|k|^r} e^{-\gamma_3'
      |k|} \ .
  \end{align*} As before,  the constant $\mathcal{D}_3'$ is
    independent of the set $\mathcal{Z}$ defining the Galerkin approximation.
\end{theorem}

Theorem \ref{thm:algebraicExponential} shows that if the initial
conditions $\omega(x,0)$ for \eqref{eq:NSvorticity} are smooth enough
then,  the solution $\omega(x,t_0)$
is real analytic for arbitrarily small time $t_0$, . Then according to theorem
\ref{thm:exponentialExponential}, it remains with in this class for
all $t>t_0$. Statements close to these were proven in the works by C.
Foias and R. Temam \cite{b:FoiasTemam89}, C. Doering and E. Titi
\cite{b:DoeringTiti95} and H. Kreiss \cite{b:Kreiss88}.  Theorem
\ref{thm:algebraicAlgebraic} is proven in \S\ref{section:pfAlgAlg}
and theorems \ref{thm:exponentialExponential} and
\ref{thm:algebraicExponential} are proven in \S
\ref{section:pfAlgExpAndExpExp}. 

The proofs of all of the theorems in this paper share a common
structure.  We consider the system of coupled ODEs for the Fourier
coefficients. Then we construct a subset $\Omega$ of the phase space
(the set of possible configurations of the Fourier modes) so that all
points in $\Omega$ possess the desired decay properties. In addition,
$\Omega$ is constructed so that it contains the initial data in its
interior.  Then we endeavor to show that the dynamics never cause
the sequence of Fourier modes to leave the subset $\Omega$.  How this
is done can be understood geometrically.  It amounts to showing that
the vector field on the boundary of $\Omega$ points into the interior
of $\Omega$. If this is true, then the solution can never escape
$\Omega$.

\section{ Proof of Theorem \ref{thm:algebraicAlgebraic} } 
\label{section:pfAlgAlg}
Fixing an arbitrary Galerkin approximation corresponding to the modes
in some finite subset $\mathcal{Z}$ of $\ZZ^2$, we write the real
version of \eqref{eq:NSFourierVorticity}. As we already mentioned, we
assume $\omega_0=0$ and, because the velocity is real, we also have
$\omega_{-k}=\bar{\omega}_k$. Setting $\omega_k=\omega^{(1)}_k + i
\omega^{(2)}_k$, we separate the equations for $\omega^{(1)}_k$ and
$\omega^{(2)}_k$ obtaining
\begin{align}
  \label{NSReal}
  \frac{d \omega_k^{(1)}}{dt}=& 2 \pi \sum_{\substack{l_1 + l_2 =k\\
    l_1, l_2 \in \mathcal{Z}}} \left[
    \omega^{(1)}_{l_1} \omega^{(2)}_{l_2} + \omega^{(2)}_{l_1}
    \omega^{(1)}_{l_2} \right] \frac{(k,l_2^\perp)}{(l_2,l_2)} - 4 \pi
    \nu |k|^\alpha \omega_k^{(1)} + g_k^{(1)}\\
 \frac{d \omega_k^{(2)}}{dt}=& -2 \pi \sum_{\substack{l_1 + l_2 =k\\
    l_1, l_2 \in \mathcal{Z}}} \left[
    \omega^{(1)}_{l_1} \omega^{(1)}_{l_2} + \omega^{(2)}_{l_1}
    \omega^{(2)}_{l_2} \right] \frac{(k,l_2^\perp)}{(l_2,l_2)} - 4 \pi
    \nu |k|^\alpha \omega_k^{(2)} + g_k^{(2)} \notag
\end{align}
where $g_k=g_k^{(1)}+ i g_k^{(2)}$.

It follows from the enstrophy estimate that $\sum_k\left[
  \bigl(w_k^{(1)}(t)\bigr)^2+ \bigl(w_k^{(2)}(t)\bigr)^2 \right] \leq
\ES^*$ and thus $\lvert w_k^{(1)}(t)\rvert \leq \sqrt{\ES^*}$ and $\lvert
w_k^{(2)}(t)\rvert \leq \sqrt{\ES^*}$ for all $k \in \ZZ^2$ and $t>0$.
Hence, for any $K_0 >0$, we can find a
$\mathcal{D}_1'=\mathcal{D}_1'(K_0)$ such that for any $t\geq 0$
$\lvert w_k^{(1)}(t)\rvert, \lvert w_k^{(2)}(t)\rvert <
\frac{\mathcal{D}_1'}{|k|^r}$ for all $k \in \ZZ^2$ with $\lvert k
\rvert \leq K_0$. We also require $\mathcal{D}_1'$ to be greater than
$\mathcal{G}$ so later estimates will arrange themselves nicely.
Recall that $\mathcal{G}(r)$ was the constant from assumption
\ref{a:powerLaw}.  Since $\mathcal{G}$ is given and only $K_0$ is ours
to vary, we will suppress the dependence of $\mathcal{D}_1'$ on
$\mathcal{G}$.

Now consider the subset 
\begin{align*}
  \Omega_1(K_0)=\left\{(\omega_k^{(1)},\omega_k^{(2)})_{ k\in \ZZ^2} :
  |\omega_k^{(j)}| \leq \frac{\mathcal{D}_1'(K_0)}{|k|^r} \mbox{ for
    all $j\in\{1,2\}$, $k \in\ZZ^2\backslash 0$} \right\} 
\end{align*}
of $(\RR^2)^{\ZZ^2}$. Its boundary is the subset 
\begin{equation*}
  \partial \Omega_{1}(K_0) = \left\{ (\omega_k^{(1)},\omega_k^{(2)})_{k\in
      \ZZ^2} : \begin{gathered}|\omega_k^{(j)}| \leq
      \frac{\mathcal{D}_1'}{|k|^r} \mbox{ for all } j \in \{1,2\},
      k \in \ZZ^2\backslash 0 \\ \mbox{ and  equality holds for some
      $\bar{k}$ and $\bar{\jmath}$ .} \end{gathered}\right\}  \ .
\end{equation*}
We shall also need the subset of this boundary 
\begin{equation*}
  \overline{\partial \Omega_{1}}(K_0) = \left\{ (\omega_k^{(1)},\omega_k^{(2)})_{k\in
      \ZZ^2} : \begin{gathered} |\omega_k^{(j)}| \leq
      \frac{\mathcal{D}_1'}{|k|^r} \mbox{ for all } j \in \{1,2\},
      k \in \ZZ^2\backslash 0 \\ \mbox{ and  equality for some
      $\bar{k}$ and $\bar{\jmath}$ with $|\bar{k}|>K_0$.} \end{gathered}\right\} \ .
\end{equation*}

Showing that the trajectories of our system remain inside of
$\Omega_1$ is equivalent to the statement of the theorem. Recall that
using the enstrophy estimate, we picked a $\mathcal{D}_1'(K_0)$ such
that if $|k| \leq K_0$ then $\lvert w_k^{(1)}(t)\rvert$ and $ \lvert
w_k^{(2)}(t)\rvert$ were bounded by $\frac{\mathcal{D}_1'}{|k|^r}$ for
all $t \in [0,\infty)$.  Thus, the only remaining way for a trajectory
to leave $\Omega_{1}(K_0)$, is through the section of the boundary
$\overline{\partial\Omega_1}(K^0)$ introduced above. Our basic idea is
to show that if $K_0$ is greater than a specific $K_{crit}$, then the
vector field on $\overline{\partial\Omega_1}(K^0)$ points inward. In
other words, the dynamics of \eqref{NSReal} can never move the system
configuration through $\partial \Omega_{1}(K_0)$.  In still
different words, $\Omega_1$ is a trapping region. Since the initial
data begins in $\Omega_1$, proving this picture would prove the
theorem.

To show that the vector field points inward, fix a point on
$\overline{\partial \Omega_{1}}(K_0)$. For definiteness, consider the case
when $\omega^{(1)}_{\bar{k}}=\frac{\mathcal{D}_1'}{|\bar{k}|^r}$ for
some $\kb$ with $|\kb| > K_0$, $|\omega_{k'}^{(1)}| \leq
\frac{\mathcal{D}_1'}{|k'|^r}$ for all $k' \in \ZZ \backslash 0$ with
$k' \neq \kb$, and $|\omega_{k}^{(2)}| \leq
\frac{\mathcal{D}_1'}{|k|^r}$ for all $k \in \ZZ^2\backslash 0$. The
other cases, namely where $\omega^{(1)}_{\bar{k}}=-
\frac{\mathcal{D}'_1}{|\bar{k}|^2}$ or $\omega^{(2)}_{\bar{k}}=\pm
\frac{\mathcal{D}'_1}{|\bar{k}|^2}$, are handled in the same manner.
We have to show that,
\begin{align}
  \label{eq:pointsIn}
  2 \pi \left\lvert \sum_{l_1 + l_2=\bar{k}} \left[\omega^{(1)}_{l_1}
      \omega^{(1)}_{l_2} + \omega^{(2)}_{l_1} \omega^{(2)}_{l_2}
    \right] \frac{(\bar{k},l_2^\perp)}{(l_2,l_2)} \right\rvert +
  \left|g_{\kb}^{(2)}\right| < 4 \pi \nu |\bar{k}|^\alpha
      \left|\omega_{\kb}^{(2)}\right| \ .
\end{align}
We shall see that the restriction that $|\bar{k}| \geq K_0 > K_{crit}$
does not depend on $\mathcal{D}_1$ only on $\ES$.

Consider the following three sums which together bound the first
abolute value on the left-hand
side of \eqref{eq:pointsIn}:
\begin{align*}
  \Sigma_1 =& \sum_{\substack{l_1+l_2=\bar{k}\\ |l_2| \leq
      |\frac{\bar{k}}{2}|} }
  \left| \left[ \omega^{(1)}_{l_1}
      \omega^{(1)}_{l_2} + \omega^{(2)}_{l_1} \omega^{(2)}_{l_2}
    \right] \right| \left|\frac{(\bar{k},l_2^\perp)}{(l_2,l_2)}\right| \\
 \Sigma_2 =& \sum_{\substack{l_1+l_2 = \bar{k}\\ |\frac{\bar{k}}{2}| <
   |l_2| \leq 2|\bar{k}| }}
  \left| \left[ \omega^{(1)}_{l_1}
      \omega^{(1)}_{l_2} + \omega^{(2)}_{l_1} \omega^{(2)}_{l_2}
    \right] \right| \left|\frac{(\bar{k},l_2^\perp)}{(l_2,l_2)}\right|
  \\
 \Sigma_3 =& \sum_{\substack{l_1+l_2=\bar{k}\\ |l_2| >
      2|\bar{k}|} }
  \left| \left[ \omega^{(1)}_{l_1}
      \omega^{(1)}_{l_2} + \omega^{(2)}_{l_1} \omega^{(2)}_{l_2}
    \right] \right| \left|\frac{(\bar{k},l_2^\perp)}{(l_2,l_2)}\right| 
\end{align*}
We treat each sum separately. For $\Sigma_1$, using the
Cauchy-Schwartz inequality
and the inequalities
${\left|\frac{(\bar{k},l_2^\perp)}{(l_2,l_2)}\right| \leq
\frac{|\bar{k}|}{|l_2|}}$, ${|l_1|\geq \frac{|\bar{k}|}{2}}$, and
$|\omega_{l_1}^{(1)}| \leq 2^r
\mathcal{D}_1'\frac{1}{|\bar{k}|^r},|\omega_{l_1}^{(2)}| \leq 2^r
\mathcal{D}_1'\frac{1}{|\bar{k}|^r}$ produces 
\begin{align}
  \label{eq:S1bound}
  \left| \Sigma_1 \right| &\leq 2^{r}
  \frac{\mathcal{D}_1'}{|\bar{k}|^r}|\bar{k}| \sum_{|l_2| \geq
    |\frac{\bar{k}}{2}|} \left[ |\omega_{l_2}^{(1)}| +
    |\omega_{l_2}^{(2)}| \right] \frac{1}{|l_2|} \notag \\
  & \leq 2^{r} \frac{\mathcal{D}_1'|\bar{k}|}{|\bar{k}|^r} \left (
    \sqrt{\sum |\omega_{l_2}^{(1)}|^2 } + \sqrt{\sum
      |\omega_{l_2}^{(2)}|^2 } \right) \sqrt{\sum_{|l_2| \leq
      |\frac{\bar{k}}{2}|} \frac{1}{|l_2|^2}} \notag \\ & \leq 2^{r+1}
  (\mbox{const}) \left(\sqrt{\ES^*}\right) |\bar{k}|
  \left(\sqrt{\ln{|\bar{k}|}}\right)
  \left(\frac{\mathcal{D}_1'}{|\bar{k}|^r}\right) \ .
\end{align}
The $(\mbox{const})$ in the final line is from the inequality
\begin{align*}
  \sum_{|l_2| \leq |\frac{\bar{k}}{2}|} \frac{1}{|l_2|^2} \leq
  (\mbox{const})^2 \ln |\bar{k}| \ .
\end{align*}

To estimate $\Sigma_2$, we use the inequalities $\left\lvert
\frac{(k,l_2^\perp)}{(l_2,l_2)}\right\rvert \leq 2$, $|\omega^{(1)}_{l_2}|
\leq 2^r \frac{\mathcal{D}_1'}{|\bar{k}|^r}$, and $|\omega_{l_2}^{(2)}|
\leq 2^r \frac{\mathcal{D}_1'}{|\bar{k}|^r}$ obtaining
\begin{align}
  \label{eq:S2bound}
  \left| \Sigma_2 \right| &\leq 2^{r+1}
  \frac{\mathcal{D}_1'}{|\bar{k}|^r} \sum_{|l_1| \leq 3 |\bar{k}|}
  \left[ |\omega_{l_1}^{(1)}| +
    |\omega_{l_1}^{(2)}| \right] \notag \\
  & \leq 2^{r+1} \frac{\mathcal{D}_1'}{|\bar{k}|^r} \left[
    \sqrt{\sum_{|l_1| \leq 3 |\bar{k}|} |\omega_{l_1}^{(1)}|^2} +
    \sqrt{\sum_{|l_1| \leq 3 |\bar{k}|} |\omega_{l_1}^{(2)}|^2}
  \right](6 |\bar{k}| +1)
  \notag \\
  & \leq 2^{r+2}  \ES  (6 |\bar{k}| +1)
  \frac{\mathcal{D}_1'}{|\bar{k}|^r} \ .
\end{align}
The factor $ (6 |\bar{k}| +1)$ arises as an estimate of the square
root of the number of lattice points $l_1 \in \ZZ^2$ for which $|l_1|
\leq 3 |\bar{k}|$.

In estimating $\Sigma_3$, we use
$\left|\frac{(\bar{k},l_2^\perp)}{(l_2,l_2)}\right| \leq
\left|\frac{\bar{k}}{l_2}\right|$ producing
\begin{align}
  \label{eq:S3bound}
  |\Sigma_3| & \leq |\bar{k}| \sum_{\substack{l_1+l_2=\bar{k} \notag \\
      |l_2|\geq 2 |\bar{k}|}} \left[ |\omega_{l_1}^{(1)}|
    |\omega_{l_2}^{(2)}| + |\omega_{l_1}^{(2)}| |\omega_{l_2}^{(1)}|
  \right] \frac{1}{|l_2|} \notag \\
  & \leq |\bar{k}| \left[ { \Biggl(\sum_{|l_1| \geq \bar{k}}
        (\omega_{l_1}^{(1)} )^2\Biggr)^\hf \Biggl( \sum_{|l_2| \geq 2
          \bar{k}} \frac{ (\omega_{l_2}^{(2)})^2}{|l_2|^2}\Biggr)^\hf} +
    {\Biggl(\sum_{|l_1|\geq \bar{k}} (\omega_{l_1}^{(2)})^2
      \Biggr)^\hf \Biggl( \sum_{|l_2| \geq
          2\bar{k}} \frac{|\omega_{l_2}^{(1)}|^2}{|l_2|} \Biggr)^\hf}\right] \notag     \\
  & \leq 2 \sqrt{\ES^*} |\bar{k}| \mathcal{D}_1' {\Bigl(\sum_{|l_2|\geq 2
      |\bar{k}|} \frac{1}{|l_2|^{2(r+1)}}\Bigr)^\hf } \leq 2 \sqrt{\ES^*}
  (\mbox{const}) |\bar{k}| \frac{\mathcal{D}_1'}{|\bar{k}|^r}
\end{align}
where $(\mbox{const})$ is defined by the inequality 
\begin{align*}
  \sum_{|l_2| \geq 2 |\bar{k}|}\frac{1}{|l_2|^{2(r+1)}} \leq
    (\mbox{const})^2 \frac{1}{|\bar{k}|^{2r}} \ .
  \end{align*}
Adding \eqref{eq:S1bound}, \eqref{eq:S2bound}, and \eqref{eq:S3bound}
together, we obtain the needed bound on the right hand side of
\eqref{eq:pointsIn}: 

\begin{align}
\label{eq:thm1FinialEst}
2 \pi \sum_{l_1+l+2=\bar{k}} |\omega_{l_1}^{(1)}||\omega_{l_2}^{(2)}|
+ | \omega_{l_1}^{(2)}| |\omega_{l_2}^{(1)}| & \leq \Big[ 2^{r+1}
(\mbox{const}) \sqrt{\ES^*} |\bar{k}| \sqrt{\ln|\bar{k}|}+ 2^{r+2} \ES^*(6
|\bar{k}|+1) \notag \\& \qquad+ 2 \sqrt{\ES^*} (\mbox{const}) |k| \Big]
\frac{\mathcal{D}_1'}{|\bar{k}|^r} \notag \\
& \leq 2^{r+2} \ES^* (\overline{\mbox{const}}) |\bar{k}|
\sqrt{\ln|\bar{k}|} \frac{\mathcal{D}_1'}{|\bar{k}|^r}
\end{align}
where $ (\overline{\mbox{const}})$ is a new constant. 

By assumption \ref{a:powerLaw} and our requirement that the
$\mathcal{D}_1'$ be greater than $\mathcal{G}$ (the constant from
assumption \ref{a:powerLaw}), we know that $|g_k|
\leq \frac{\mathcal{D}_1'}{|k|^{r-\alpha+\epsilon}}$.  Thus, inequality
\eqref{eq:pointsIn} will be satisfied if
\begin{align}
  \label{eq:veryLastThm1}
  \left[  2^{r+2} \ES^*  (\overline{\mbox{const}}) \frac{|\bar{k}|
    \sqrt{\ln|\bar{k}|}}{|\kb|^\alpha}  + \frac{1}{|\kb|^\epsilon}
    \right]\frac{\mathcal{D}_1'}{|\bar{k}|^{r-\alpha}} \leq  4 \pi \nu
    \frac{\mathcal{D}_1'}{|\kb|^{r-\alpha}} \ .
\end{align}
From this we see that for all $\alpha>1$, there exists $K_{crit}$ so that if $|\bar{k}|
\geq K_{crit}$ then \eqref{eq:veryLastThm1} holds. Also notice that
$K_{crit}$ is independent of our choice of $\mathcal{D}_1'$
except for the condition that $\mathcal{D}_1' > \mathcal{G}$.
Thus we can find $K_{crit}$ first and then fix $K_0$ which determines
$\mathcal{D}_1'$.

\section{ Proofs of Theorems \ref{thm:exponentialExponential} and \ref{thm:algebraicExponential} }
\label{section:pfAlgExpAndExpExp}

We begin by stating the central estimate on which both theorems rely.  It
requires estimates similar in spirit to the previous theorem and will be
proven at the end of the section. We present a $d$-dimensional version of
the lemma because it will be useful in the discussions of the 3-dimensional
setting in the next section.
\begin{lemma}
\label{l:sum}  Let $\{a_k\}$ and $\{b_k\}$ be two sequences with $k \in \ZZ^d$.  If
  for some $r >d-1$ and some $\mathcal{C}>0$
\begin{align*}
  |a_{k}| \leq \frac{\mathcal{C}}{|k|^r} \qquad |b_k| \leq \frac{\mathcal{C}}{|k|^r}
\end{align*}
then for all $k \in \ZZ^d$
  \begin{align*}
    \sum_{\substack{l_1+l_2=k\\ l_1, l_2 \in \ZZ^d}} \left\lvert a_{l_1}\right\rvert \left\lvert
      b_{l_2}\right\rvert \frac{|k|}{|l_2|} \leq (\mbox{const})
    \left( 2^{r} |k| + 2^{r+1}(6|k|+1)^{\frac{d}{2}} + \frac{1}{2} |k|^{d-1-r} \right)
       \frac{\mathcal{C}^2}{|k|^r}
  \end{align*}
where the constant depends only on $r$ and not on $k$.
\end{lemma}
We now turn to the proof of theorem \ref{thm:exponentialExponential}.
\begin{proof}[Proof of theorem \ref{thm:exponentialExponential}]
If $|\omega_k^{(1)}(0)| \leq \frac{\mathcal{D}_2}{|k|^r}
e^{-\gamma_2|k|}$, $|\omega_k^{(2)}(0)|\leq
\frac{\mathcal{D}_2}{|k|^r} e^{-\gamma_2|k|}$ then surely
$|\omega_k^{(1)}|, |\omega_{k}^{(2)}| \leq
\frac{\mathcal{D}_2}{|k|^r}$. Therefore by theorem
\ref{thm:algebraicAlgebraic}, one can find a constant
$\bar{\mathcal{D}}_2$ such that $|\omega_k^{(1)}(t)|,
|\omega_k^{(2)}(t)| \leq \frac{\bar{\mathcal{D}}_2}{|k|^r}$ for all $k
\in \ZZ^2 \backslash \{ 0 \}$. Let us set $\mathcal{D}'_2 = \max(2
\bar{\mathcal{D}}_2,\mathcal{G})$ where $\mathcal{G}$ is the constant
from assumption \ref{a:exponential}. The numerical factor 2 is somewhat
arbitrary. We could chose any factor greater than 1; we take 2 for
simplicity.

Choose $K_0 \geq 0$ and consider the set
\begin{align*}
  \Omega_2(K_0) = \left\{ (\omega_k^{(1)}, \omega_k^{(2)})_{k\in \ZZ^2
  \backslash\{0\}} :  \begin{gathered}|\omega_k^{(1)}| \leq
  \frac{\mathcal{D}_2'}{|k|^r}e^{-\gamma_2'|k|},  |\omega_k^{(2)}|
  \leq\frac{\mathcal{D}_2'}{|k|^r}e^{-\gamma_2'|k|} 
\end{gathered} 
\right\}
\end{align*}
The value of $\gamma_2'=\gamma_2'(K_0)$ is chosen in such a way that
the inequalities $|\omega_k^{(i)}(t)| \leq
\frac{\bar{\mathcal{D}}_2}{|k|^r}$ given by theorem
\ref{thm:algebraicAlgebraic} imply that $|\omega_k^{(i)}(t)|\leq
\frac{\mathcal{D}_2'}{|k|^r} e^{-\gamma_2' |k|}$ for all $k$, $|k| \leq
K_0$, and that for $|k| \geq K_0$, $e^{-\gamma_2'|k|}\geq
e^{-\gamma|k|^{1+\delta}}$. Here $\gamma$ and $\delta$ are the
constants from assumption \ref{a:exponential}.

As in the proof of theorem \ref{thm:algebraicAlgebraic}, we shall show
that for sufficiently large $K_0$ the vector field corresponding to
\eqref{eq:pointsIn} is directed inside $\Omega_2(K_0)$ along the part
of the boundary $\partial \Omega_2(K_0)$ where $|\omega_k^{(i)}| \leq
\frac{\mathcal{D}'_2}{|k|^r} e^{-\gamma_2'|k|}$ for all $k\in
\ZZ^2\backslash\{0\}$ with $|k|\geq K_0$ and for at least one of
these, say
$\bar{k}$,  we have  equality. It will
be shown that our restriction from below on $K_0$, needed to ensure
the vector field points inward, will not depend on $\gamma_2'$. This
will yield the stated result.

As in theorem \ref{thm:algebraicAlgebraic}, consider for definiteness
the case where
$\omega^{(1)}_{\bar{k}}=\frac{\mathcal{D}'_2}{|\bar{k}|^2}
e^{-\gamma_2'|\bar{k}|}$. The other cases are handled in the same
manner. As before, we have to show that the vector field points
inward. This would be assured if
\begin{align}
  \label{eq:inWard2}
  2\pi \left\lvert \sum_{l_1+l_2=\bar{k}} [ \omega_{l_1}^{(1)}
  \omega_{l_2}^{(2)} + \omega_{l_2}^{(2)} \omega_{l_2}^{(1)} ]
  \frac{(\kb,l^\perp_2)}{(l_2,l_2)}\right\rvert +|g_{\kb}| < 4 \pi^2 |\kb|^\alpha
  \frac{\mathcal{D}_2'}{|\kb|^r} e^{-\gamma_2'|\kb|} \ .
\end{align}
This time we do not use the enstrophy estimate as previously. Instead, we use the
estimates
$|\omega^{(1)}_{l_1}| \leq \frac{\mathcal{D}_2'}{|l_1|^r}
e^{-\gamma_2'|l_1|}$ and $|\omega^{(2)}_{l_1}| \leq \frac{\mathcal{D}_2'}{|l_2|^r}
e^{-\gamma_2'[|l_2|}$ .

Let us put $v_k^{(j)}= e^{\gamma_2'|k|} \omega_k^{(j)}$, $j=1,2$, $k\in
\ZZ^2 \backslash 0$. In terms of $v$, \eqref{eq:inWard2} becomes
\begin{align}
  \label{inWard3}
  2\pi \left\lvert \sum_{l_1+l_2=\bar{k}} [ v_{l_1}^{(1)}
    v_{l_2}^{(2)} + v_{l_2}^{(2)} v_{l_2}^{(1)} ]
    \frac{e^{-\gamma_2'|l_1| - \gamma_2'|l_2|}}{e^{-\gamma_2'|\kb|}}
    \frac{(\kb,l^\perp_2)}{(l_2,l_2)}\right\rvert +
  |g_{\kb}|e^{\gamma_2'|\kb|} < 4 \pi^2 |\kb|^\alpha
  \frac{\mathcal{D}_2'}{|\kb|^r} \ .
\end{align}
First notice that $\frac{e^{-\gamma_2'|l_1|}e^{-
    \gamma_2'|l_2|}}{e^{-\gamma_2'|\kb|}} \leq 1$ so it may be neglected.
Second notice that for $v_k^{(j)}$, we have the estimate $|v_k^{(j)}|
\leq \frac{\mathcal{D}_2'}{|k|^r}$ for $k \in \ZZ^2 \backslash 0$.
Lastly, we know that $\lvert \frac{(\kb,l^\perp_2)}{(l_2,l_2)} \rvert
\leq \frac{|k|}{|l_2|}$. These estimates allow us to apply lemma \ref{l:sum}, producing
\begin{align}
  \label{eq:inWard4}
  2 \pi \Big\lvert \sum \left[v_{l_1}^{(1)} v_{l_2}^{(2)} +
      v_{l_2}^{(2)} v_{l_2}^{(1)} \right] &    \frac{e^{-\gamma_2'|l_1|
      - \gamma_2'|l_2|}}{e^{-\gamma_2'|\kb|}}
      \frac{(\kb,l^\perp_2)}{(l_2,l_2)} \Big\rvert \\ & \leq 2 \pi \ 
  \mbox{const} \ \left( 2^{r+1} |\kb| + 2^{r+2}(6|\bar{k}|+1) + |\kb|^{1-r}
  \right) \mathcal{D}_2' \frac{\mathcal{D}_2'}{|\kb|^r} \ .\notag 
\end{align}
From this estimate, we see that if
\begin{align}
  \label{eq:inWard5}
  2 \pi (\mbox{const})\left( 2^{r+1} |\kb| + 2^{r+2}(6|\bar{k}|+1) + 2 |\kb|^{1-r}
  \right) \mathcal{D}_2' + \frac{\mathcal{G}}{\mathcal{D}_2'}
  \frac{e^{-\gamma|\kb|^{1+\delta}}}{e^{-\gamma_2'|\kb|}}|\kb|^{\alpha-\epsilon}
  < 4 \pi^2 \nu |\kb|^\alpha
\end{align}
then we have established \eqref{inWard3}, which was our goal. Notice
that we chose $\mathcal{D}_2' \geq \mathcal{G}$ and $\gamma_2'$ such
that $\frac{e^{-\gamma|\kb|^{1+\delta}}}{e^{-\gamma_2'|\kb|}} \leq 1 $
for all $k$ with $|k| \geq K_0$. Since $\alpha>1$ by picking $K_0$ large enough, we can
force \eqref{eq:inWard5} to hold. This is the criteria which sets the level of
$K_{crit}$.  The proof of theorem \ref{thm:exponentialExponential} is concluded.
\end{proof}
We now present the proof of theorem
\ref{thm:algebraicExponential}. Its structure is very similar to the
previous proof and also employs lemma \ref{l:sum} but uses a slightly
different change of variable.
\begin{proof}[Proof of theorem \ref{thm:algebraicExponential}]
  Let $\mathcal{D}_1'$ be the constant given by theorem
  \ref{thm:algebraicAlgebraic}, that is such that $|\omega_k(t)| \leq
  \frac{\mathcal{D}_1'}{|k|^2}$ for all $k \in \ZZ^2 \backslash 0$ and
  all $t$. Let us put $v_k^{(j)}=\omega_k^{(j)} e^{\gamma_3 t |k|},
  j=1,2$ where the constant $\gamma_3$ will be determined later. The
  evolution of the $v_k^{(1)}(t)$ are described by the following ODEs
\begin{align}
  \label{eq:vkSystem}
  \frac{d v_k^{(1)}(t)}{dt}= & \gamma_3 |k| v_k^{(1)}(t) - 4 \pi^2 \nu
  |k|^\alpha v_k^{(1)}(t) +g_k^{(1)} e^{\gamma_3 t |k|}\\ &- 2 \pi
  \sum_{l_1+l_2=k} \left[ v_{l_1}^{(1)}(t) v_{l_2}^{(2)}(t)+
    v_{l_1}^{(2)}(t)v_{l_2}^{(1)}(t) \right] \frac{e^{-\gamma_3 t
      |l_1|}e^{-\gamma_3 t |l_2|}}{e^{-\gamma_3 t |k|}} 
  \frac{(k,l_2^\perp)}{(l_2,l_2)}  \ . \notag
\end{align}
The analogous equations describe the evolution of the $v_k^{(2)}(t)$.

The methods of the previous section can be applied to this coupled
system. We fix a time $t_0>0$ and an arbitrary positive constant
$\gamma_0$. For $t=0$, we have the inequalities
\begin{align*}
  |v_k^{(1)}(0)| \leq \frac{\mathcal{D}_3}{|k|^r} \qquad
   |v_k^{(2)}(0)| \leq \frac{\mathcal{D}_3}{|k|^r}
\end{align*}
for all $k$. In light of the definition of $v_k(t)$, theorem
\ref{thm:algebraicExponential} would be proven if we show that 
\begin{align}\label{eq:endT3}
  |v_k^{(1)}(t_0)| \leq \frac{\mathcal{D}_3'}{|k|^r} \qquad
   |v_k^{(2)}(t_0)| \leq \frac{\mathcal{D}_3'}{|k|^r}
\end{align}
for some appropriate $\mathcal{D}_3'$.

As in the proof of theorem \ref{thm:algebraicExponential}, we put
$\mathcal{D}_3'=\max(2\mathcal{D}_1',\mathcal{G})$ where
$\mathcal{G}$ is again the constant from assumption \ref{a:exponential}. For any fixed
$K_0$, we can find a $\gamma_3$ so that the following three conditions
hold. First, the inequalities $|\omega_k^{(j)}(t)| \leq
\frac{\mathcal{D}_1'}{|k|^r}$, imply $|v_k^{(j)}(t)| \leq
\frac{\mathcal{D}_3'}{|k|^r}$ for $j=1,2$, $t\in [0,t_0]$, and $|k|
\leq K_0$. Second, so $e^{\gamma_3 t |k|} \leq e^{\gamma
  |k|^{1+\delta}}$ for $k\in \ZZ^2$ with $|k| \geq K_0$ and $t \in [0,t_0]$. In
this condition the constants $\gamma$ and $\delta$ are again from
assumption $\ref{a:exponential}$. Third, we can always assume that
$\gamma_3 \leq \gamma_0$. (This last assumption is to simplify the
exposition and is not really needed as $\gamma_3$ decreases as we
increase $K_0$.)

Now consider the set 
\begin{align*}
  \Omega_3(K_0)= \left\{ \left(v_k^{(1)},v_k^{(2)}\right)_{k\in
  \ZZ^2\backslash 0} \mbox{ with } |v_k^{(j)}| \leq
  \frac{\mathcal{D}_3'}{|k|^r} \mbox{ for } j=1,2 \mbox{ and } |k|
  > K_0 \right\} \ .
\end{align*}

Again we will show that if $K_0$ is greater than some $K_{crit}$, the
vector field along the boundary of $\Omega_3(K_0)$ points inward.  The
calculation parallels that in theorem
\ref{thm:exponentialExponential}. For definiteness, we assume that
$v_{\kb}^{(1)}(t)=\frac{\mathcal{D}_3'}{|\kb|^r}$ for some $\kb$ with $|\kb|
> K_0$ and that the inequality bounds which define $\Omega_3$ hold for
all other $k$. The other cases proceed analogously.

We wish to show that the vector field points inward. Since $\gamma_3
\leq \gamma_0$, from \eqref{eq:vkSystem}, we see that it is sufficient
to show that for $t \in [0,t_0]$
\begin{align}
\label{T3in}
  (4 \pi^2 \nu |\kb|^\alpha- \gamma_0|\kb|) v^{(1)}_{\kb} &>
  2 \pi \left| \sum_{l_1+l_2=\kb} \left[ v_{l_1}^{(1)}(t)
    v_{l_2}^{(2)}(t)+ v_{l_1}^{(2)}(t)v_{l_2}^{(1)}(t) \right]
    \frac{(k,l_2^\perp)}{(l_2,l_2)} \right| \notag \\ & \qquad +
  |g_{\kb}^{(1)}|e^{\gamma_3 t |\kb|}  \ .
\end{align}
Here, as before, we have neglected the factor $\frac{e^{-\gamma_3 t
    |l_1|}e^{-\gamma_3 t |l_2|}}{e^{-\gamma_3 t |k|}}$ as it is always
less than 1. After applying the inequalities $\mathcal{G} \leq
\mathcal{D}_3'$, $e^{\gamma_3 t |k|} \leq e^{\gamma |k|^{1+\delta}}$
  and lemma \ref{l:sum}, we see that \eqref{T3in} holds if
\begin{align}
  \label{T3in2}
  4 \pi^2 \nu > \gamma_0 \frac{|\kb|}{|\kb|^\alpha} +
  (\mbox{const})\mathcal{D}_3'\left[ 2^{r+1}  \frac{|\kb|}{|\kb|^\alpha}
    + 2^{r+2} \frac{7|\kb|}{|\kb|^\alpha} +
    \frac{|\kb|^{1-r}}{|\kb|^\alpha}\right] +
  \frac{\mathcal{G}}{\mathcal{D}_3'}  \frac{1}{|\kb|^\alpha} \ .
\end{align}
Because $\alpha >1$ and $r>2$, by making $\kb$ large enough we can
force \eqref{T3in2} to hold. This shows that the solution to any
Galerkin approximation stays in $\Omega_3$ until the time $t_0$ and
thus \eqref{eq:endT3} holds and the proof is complete.
\end{proof}

\begin{proof}[Proof of Lemma \ref{l:sum}:] 
  As in the proof of theorem \ref{thm:algebraicAlgebraic}, we estimate
  separately three sums. 
\begin{align*}
  \Sigma_1=&\sum_{\substack{|l_2| \leq |\frac{k}{2}| \\
      l_1+l_2=k}} \left\lvert a_{l_1}\right\rvert \left\lvert b_{l_2}
  \right\rvert 
    \frac{|k|}{|l_2|} \\
  \Sigma_2=&\sum_{\substack{|\frac{k}{2}| < |l_2| \leq 2|k| \\
      l_1+l_2=k}} \left\lvert a_{l_1} \right\rvert\left\lvert b_{l_2}
  \right\rvert
   \frac{|k|}{|l_2|}  \\
  \Sigma_3=&\sum_{\substack{|l_2| > 2|k| \\ l_1+l_2=k}} \left\lvert
    a_{l_1}\right\rvert \left\lvert b_{l_2}\right\rvert \frac{|k|}{|l_2|}
\end{align*}

Since in $\Sigma_1$, the norm of $|l_1| \geq |\frac{\kb}{2}|$, we can
write
\begin{align}
  \label{eq:T2S1}
  \Sigma_1 \leq \sum_{|l_2| \leq |\frac{k}{2}|} |a_{l_1}| |b_{l_2}|
  \frac{|k|}{|l_2|} \leq \frac{2^{r} \left( \mathcal{C}
    \right)^2}{|k|^r} |k| \sum_{{|l_2| \leq
      |\frac{k}{2}|}}\frac{1}{|l_2|^{r+1}} \leq 2^{r} (\mbox{const})
  |k| \frac{\mathcal{C}^2}{|k|^r}
\end{align}
where the constant is defined by the inequality 
\begin{align*}
  \sum_{l_2 \in \ZZ^d \backslash 0 } \frac{1}{|l_2|^{r+1}} \leq \mbox{const} \ .
\end{align*}
For this sum to be finite, we need $r+1> d$. For $\Sigma_2$ we have
$\frac{|k|}{|l_2|} \leq 2$ and hence
\begin{align}
  \label{eq:T2S2}
  \Sigma_2 & \leq 2 \sum_{|\frac{k}{2}| < |l_2| \leq 2|k|}
  \frac{\mathcal{C}^2}{|l_1|^r |l_2|^r} \leq \frac{2^{r+2} 
    (\mathcal{C})^2}{|\kb|^r} \sum_{|l_1|\leq 3|k|} \frac{1}{|l_1|^r} \notag
  \\ & \leq  \frac{2^{r+2} 
    (\mathcal{C})^2}{|\kb|^r} \Big(\sum_{|l_1|\leq 3|k|}
  \frac{1}{|l_1|^{2r}} \Big)^\hf \Big(\sum_{|l_1|\leq 3|k|} 1  \Big)^\hf \leq
  2^{r+1} (\mbox{const})\frac{\mathcal{C}^2(6|k|+1)}{|k|^r} \ .
\end{align}
Here the constant is the absolute constant defined by 
\begin{align*}
  \sum_{l_1 \in \ZZ^d\backslash 0} \frac{1}{|l_1|^{2r}} \leq  \mbox{const} \ .
\end{align*}
For this sum to be finite, we need $2r> d$.  For $\Sigma_3$ we have
$\frac{|k|}{|l_2|} \leq
\frac{1}{2}$.  Hence, we can write
\begin{align}
  \label{eq:T2S3}
  \Sigma_3 &\leq \frac{1}{2} \sum_{ |l_2| \geq 2|\kb| }
  \frac{\mathcal{C}^2}{|l_1|^r |l_2|^{r+1}} \leq
  \frac{\mathcal{C}^2}{2} \Biggl(\sum_{\substack{|l_1| > |k|\\l_1 \in
  \ZZ^d\backslash 0}}
    \frac{1}{|l_1|^{2r}} \Biggr)^\hf \Biggl(\sum_{\substack{|l_2| \geq 2|k|\\l_2 \in
  \ZZ^d\backslash 0}}
    \frac{1}{|l_2|^{2r+2}}\Biggr)^\hf
  \notag \\
  & \leq
  (\mbox{const})\frac{|k|^{d-1-r}}{2}\frac{\mathcal{C}^2}{|k|^r} \ .
\end{align}

Collecting together \eqref{eq:T2S1},\eqref{eq:T2S2},\eqref{eq:T2S3},
we obtain the lemma.
\end{proof}

\section{The three-dimensional setting}

This paper is mainly concerned with presenting an elementary proof
of existence and uniqueness results in the two-dimensional setting. However,
these techniques can also be used to gain some insight into the
three-dimensional setting. On the three torus, the Navier-Stokes
equations take the form
\begin{align}
  \label{eq:3DNS}
  \frac{\partial u_i}{\partial t} + \sum_{j=1,2,3} u_j
  \frac{\partial u_i}{\partial x_j} &= \nu \Delta u_j
  -\frac{\partial p}{\partial x_i} + f_i \qquad i=1,2,3
  \notag \\
  \sum_{i=1,2,3}  \frac{\partial u_i}{\partial x_i} &= 0
\end{align}
where $\nu > 0$ is again the viscosity, $p$ is the pressure, and the $f_i$ are the
components of the external, time-dependent forcing. As before, we introduce the
vorticity $\omega(x,t)=(\omega_1(x,t),\omega_2(x,t),\omega_3(x,t))=
(\frac{\partial u_2}{\partial x_3} -\frac{\partial u_3}{\partial x_2},
\frac{\partial u_3}{\partial x_1} -\frac{\partial u_1}{\partial
  x_3},\frac{\partial u_1}{\partial x_2} -\frac{\partial u_2}{\partial
  x_1})$. The vorticity obeys the equation 
\begin{align}
  \label{3DvortEqn}
  \frac{\partial \omega_i}{\partial t} + \sum_j u_j \frac{\partial
  \omega_i}{\partial x_j} =  \sum_j \omega_j \frac{\partial
  u_i}{\partial x_j}+ \nu \Delta \omega_i + g_i \qquad i=1,2,3
\end{align}
where the $g_i$ are the components of $\mbox{curl}f$. Moving to
Fourier space where 
\begin{align*}
  u(x,t)=\sum_{k \in \ZZ^3} u_k(t) e^{2 \pi i (k,x)} \ \mbox{ and } \ 
  \omega(x,t)=\sum_{k \in \ZZ^3} \omega_k(t) e^{2 \pi i (k,x)} \ ,
\end{align*}
we obtain
\begin{multline}
  \label{3DvorticityFourier}
  \frac{d \omega_k(t)}{dt} = - 2 \pi i \sum_{l_1+l_2=k} \Bigl[
    \bigl(u_{l_1}(t),l_2\bigr)\omega_{l_2}(t) -
    \bigl(\omega_{l_1}(t),l_2\bigr)u_{l_2}(t) \Bigr] - 4 \pi^2 \nu
    |k|^2 \omega_k(t) + g_k(t) \ .
\end{multline}
Here the $g_k(t)$ are the Fourier components of the forcing
$g(x,t)$. In addition, we can replace the Laplacian with
the more general differential operator $|\nabla|^\alpha$ with $\alpha>1$.

The incompressibility condition implies that
\begin{align}
  \label{incompress}
  u_k(t) \perp k
\end{align}
for all $k \in \ZZ^3$. Similarly, it follows that $ \omega_k(t) \perp
k$, $\omega_k(t) \perp u_k(t)$, and $|\omega_k(t)|=|k| |u_k(t)|$.
Hence, $(k,u_k,\omega_k)$ is a right-handed orthogonal (but not
orthonormal) frame.

Since $(u_{l_1}(t),l_1)=(\omega_{l_1}(t),l_1)=0$, we can rewrite
\eqref{3DvorticityFourier} as
\begin{multline}
  \label{3DvorticityFourier2}
  \frac{d \omega_k(t)}{dt} = - 2 \pi i \sum_{l_1+l_2=k} \Bigl[
    \bigl(u_{l_1}(t),k\bigr)\omega_{l_2}(t) -
    \bigl(\omega_{l_1}(t),k\bigr)u_{l_2}(t) \Bigr]  - 4 \pi^2 \nu
    |k|^\alpha \omega_k(t) + g_k(t) \ .
\end{multline}
As before, we begin by restricting our attention to a finite subset
$\mathcal{Z} \subset \ZZ^3$. The finite-dimensional Galerkin system
corresponding to $\mathcal{Z}$ is
\begin{multline}
  \label{3DvorticityFourier3}
  \frac{d \omega_k(t)}{dt} = - 2 \pi i \sum_{\substack{l_1+l_2=k\\
    l_1,l_2 \in \mathcal{Z}}} \Bigl[
    \bigl(u_{l_1}(t),k\bigr)\omega_{l_2}(t) -
    \bigl(\omega_{l_1}(t),k\bigr)u_{l_2}(t) \Bigr] - 4 \pi^2 \nu
    |k|^\alpha \omega_k(t) + g_k(t) \ . \tag{\problemB}
\end{multline}
Furthermore, to simplify the arguments, we assume that the forcing
$g(x,t)$ is a trigonometric polynomial which implies that all but a
finite number of the $g_k$ are identically zero. We will always
analyze wave numbers above the band which is directly forced; hence,
we may neglect the $g_k$. This is only for convenience. The forcing
can be included in the same way as it was in the two-dimensional
setting.

Our development is based upon the basic energy estimate (see
\cite{b:CoFo88,b:DoGi95,b:Temam79}). It states that given any initial data such
that $\sum_{k \in \ZZ^3} |u_k(0)|^2 = E_0 < \infty$ then there exists
a constant $E^*$ depending only on $E_0$, $\nu$, $\sup_t
|g(\cdot,t)|_{L^2}$ such that for any finite-dimensional Galerkin
approximation, defined by $\mathcal{Z} \subset \ZZ^3$, we have $\sum_{k
  \in \mathcal{Z}} |u_k(t)|^2 < E^*$ for all $t>0$.

When $\alpha=2$, the system (\problemB) corresponds to the
Navier-Stokes equations. Unfortunately, we are unable to prove the
theorems in this setting analogous to theorems \ref{thm:algebraicAlgebraic},
\ref{thm:exponentialExponential}, and \ref{thm:algebraicExponential}
when  $\alpha=2$. However, if we increase $\alpha$, we can.

\begin{theorem}
 \label{thm3D1}
 Consider the system (\problemB) with an $\alpha>2.5$ and satisfying
 assumption \ref{a:powerLaw}. If the initial data $\{\omega_k(0)\}$
 are such that
  \begin{align*}
    |\omega_k(0)| \leq \frac{\mathcal{D}_4}{|k|^r}
  \end{align*}
  for all $k\in\ZZ^3$ with $r>1.5$ then there exists a constant
  $\mathcal{D}_4'$, independent of $\mathcal{Z}$, so that 
  \begin{align*}
    |\omega_k(t)| \leq \frac{\mathcal{D}_4'}{|k|^r} 
  \end{align*}
 for all $k\in\ZZ^3$ and $t\geq 0$. 
\end{theorem}

\begin{theorem}
 \label{thm3D2}
  Consider the system (\problemB) with an $\alpha>2.5$ and satisfying
  assumption \ref{a:exponential}. If the initial
  data $\{\omega_k(0)\}$ are such that
  \begin{align*}
    |\omega_k(0)| \leq \frac{\mathcal{D}_5}{|k|^r} e^{-\gamma_5|k|}
  \end{align*}
  for all $k\in\ZZ^3$ with $r>2$ then there exists  constants
  $\mathcal{D}_5'< \infty$ and $\gamma_5'>0$, both independent of $\mathcal{Z}$, so that 
  \begin{align*}
    |\omega_k(t)| \leq \frac{\mathcal{D}_5'}{|k|^r}  e^{-\gamma_5'|k|}
  \end{align*}
 for all $k\in\ZZ^3$ and $t\geq 0$. 
\end{theorem}

\begin{theorem}
 \label{thm3D3}
 Consider the system (\problemB) with an $\alpha>2.5$ and satisfying
 assumption \ref{a:exponential}. If the initial data $\{\omega_k(0)\}$
 are such that
  \begin{align*}
    |\omega_k(0)| \leq \frac{\mathcal{D}_6}{|k|^r} 
  \end{align*}
  for all $k\in\ZZ^3$ with $r>2$ then for any $t_0>0$ there exists  constants
  $\mathcal{D}_6'< \infty$ and $\gamma_6'>0$, both independent of $\mathcal{Z}$, so that 
  \begin{align*}
    |\omega_k(t_0)| \leq \frac{\mathcal{D}_6'}{|k|^r}  e^{-\gamma_6'|k|}
  \end{align*}
 for all $k\in\ZZ^3$.
\end{theorem}

Of these three theorems, we will only give the proof of the first. The
second two will be the consequence of two more general theorems given
below. They apply to all $\alpha >1.5$ but require the additional
assumption that the enstrophy of all Galerkin approximations, starting
from a given set of initial data, remains uniformly bounded in time. This
is not known in general. However, when $\alpha>2.5$, theorem
\ref{thm3D1} implies this bound. Hence, the two theorems below apply
to (\problemB) when $\alpha>2.5$ without any assumption on $\ES(t)$. In
light of theorem \ref{thm3D1}, theorem \ref{3DifOnlyItWereTrue} and
\ref{3DifOnlyItWereTrue2} respectively yield theorem \ref{thm3D2}
and  \ref{thm3D3}  when  $\alpha>2.5$.

\begin{theorem}
\label{3DifOnlyItWereTrue}
Let $\{u_k(t)\}$ be a solution to (\problemB) with $\alpha>1.5$ such
that $\sum_{\ZZ^3} |\omega_k(t)|^2 < \ES^* < \infty$ for all $t>0$. If
$|\omega_k(0)| \leq \frac{\mathcal{D}_7}{|k|^r}$ for some
$\mathcal{D}_7 <\infty$ and $r>2$ then for any $t_1 >0$ there exists a
$\gamma_7>0$ and $\mathcal{D}_7'< \infty$ such that
  \begin{align*}
    |\omega_k(t_1)| \leq \frac{\mathcal{D}_7'}{|k|^r}e^{-\gamma_7|k|} \ .
  \end{align*}
\end{theorem}

\begin{theorem}
\label{3DifOnlyItWereTrue2}
  Let $\{u_k(t)\}$ be a solution to (\problemB) with $\alpha>1.5$ such
  that $\sum_{\ZZ^3} |\omega_k(t)|^2 < \ES^* < \infty$ for all $t>0$. If for
  some $\mathcal{D}_8 <\infty$, $\gamma_8>0$, and $r>2$, $|\omega_k(0)| \leq
  \frac{\mathcal{D}_8}{|k|^r}e^{-\gamma_8 |k|}$  then there exists a
  $\gamma_8'>0$ and $\mathcal{D}_8'< \infty$ such that for all $t>0$
  \begin{align*}
    |\omega_k(t)| \leq \frac{\mathcal{D}_8'}{|k|^r}e^{-\gamma_8'|k|} \ .
  \end{align*}

\end{theorem}
The above two theorems apply to (\problemB) for $\alpha>1.5$. In
particular, this means that they cover the standard Navier-Stokes
equation which corresponds to $\alpha=2$.  (One can lower the
restriction on $\alpha$ to $\alpha>1$ at the cost of raising the
restriction on $r$ to $r>3$. Similarly, one lowers the restriction on
$r$ to $r>1.5$ at the cost of making $\alpha>2.5$.)

In proving these two theorems, it was necessary to assume that
$\sum_{\ZZ^3} |\omega_k(t)|^2$ remained uniformly bounded in time.
Without such an assumption, we are forced to consider only $\alpha$
which do not correspond to the Navier-Stokes equation. Notice that
theorem \ref{thm3D1} implies that $\sum_{k\in \ZZ^3} |\omega_k(t)|^2 <
\mbox{const} < \infty$ for all $t>0$ and hence theorems
\ref{3DifOnlyItWereTrue} and \ref{3DifOnlyItWereTrue2} apply showing
that the solution is analytic after $t=0$.

In proving the above results, it is again convenient to split the system
(\problemB) into the equations for the real and imaginary parts of
$\{u_k\}_k$ and $\{\omega_k\}_k$.  Letting $u_k(t)=u^{(1)}_k(t)+ i
u^{(2)}_k(t)$, $\omega_k(t)=\omega_k^{(1)}(t)+i \omega_k^{(2)}(t)$,
and $g_k(t)=g^{(1)}_k(t)+ i g^{(2)}_k(t)$; we obtain
\begin{align}
  \label{3DVorticityRealImag}
  \frac{d \omega_k^{(1)}(t)}{dt} = &2 \pi \sum_{\substack{l_1+l_2=k\\
      l_1 , l_2 \in \mathcal{Z}}} \Bigl[ \bigl(
  u^{(1)}_{l_1}(t),k\bigr)\omega^{(2)}_{l_2}(t)+
  \bigl(u^{(2)}_{l_1}(t),k\bigr)\omega^{(1)}_{l_2}(t) - \bigl(
  \omega^{(2)}_{l_1}(t),k\bigr)u^{(1)}_{l_2}(t)  -
  \bigl( \omega^{(1)}_{l_1}(t),k\bigr)u^{(2)}_{l_2}(t) \Bigr]
  \notag \\ & \qquad -4 \pi^2 \nu |k|^\alpha \omega_k^{(1)}(t) + g_k^{(1)}(t) \tag{\problemBr}\\
  \frac{d \omega_k^{(2)}(t)}{dt} = &-2 \pi \sum_{\substack{l_1+l_2=k\\ 
      l_1 , l_2 \in \mathcal{Z}}} \Bigl[ \bigl(
  u^{(1)}_{l_1}(t),k\bigr)\omega^{(1)}_{l_2}(t)-
  \bigl(u^{(2)}_{l_1}(t),k\bigr)\omega^{(2)}_{l_2}(t) - \bigl(
  \omega^{(1)}_{l_1}(t),k\bigr)u^{(1)}_{l_2}(t)  +
  \bigl( \omega^{(2)}_{l_1}(t),k\bigr)u^{(2)}_{l_2}(t) \Bigr] \notag
      \\ & \qquad + g_k^{(2)}(t)  -4 \pi^2
  \nu |k|^\alpha \omega_k^{(2)}(t) \tag{\problemBi}
\end{align}

\begin{proof}[Proof of Theorem \ref{thm3D1}]
  By the energy estimate, we know that $|u_k^{(j)}(t)| \leq
  \sqrt{E^*}$ for all $t\geq 0$ and $j=1,2$. Hence,
  $|\omega^{(j)}_k(t)| \leq |k| \sqrt{E^*}$. Fixing a $K_0$, set
  $\mathcal{D}_4'(K_0)=K_0 \mathcal{D}_4$. With this choice,
  $|\omega_k^{(j)}(t)| \leq \mathcal{D}_4'(K_0)$ for all $t \geq 0$,
  $j=1,2$, and $k \in \ZZ^3$ with $|k| \leq K_0$.  As before, consider
  the set
  \begin{align*}
     \Omega_4(K_0)=\left\{(\omega_k^{(1)},\omega_k^{(2)})_{ k\in \ZZ^3} :
  |\omega_k^{(1)}| \leq \frac{\mathcal{D}_4'(K_0)}{|k|^r},|\omega_k^{(2)}|
  \leq \frac{\mathcal{D}_4'(K_0)}{|k|^r} \mbox{ for all $k$, $|k| > K_0$}
\right\} \ .
  \end{align*}
We have to show that if $K_0$ is taken to be sufficiently
large, the vector field points inward along $\partial \Omega_4$.

We pick a point on the boundary. For definiteness, we will again consider
the case when $\omega_{\kb}^{(1)}=\frac{\mathcal{D}_4'}{|\kb|^r}$
and $\omega_{\kb}^{(2)}\leq\frac{\mathcal{D}_4'}{|\kb|^r}$
for some $\kb$ with $|\kb|\geq K_0$ and
$\omega_{k}^{(j)}\leq\frac{\mathcal{D}_4'}{|k|^r}$ for all other $k$
with $k \neq \kb$.
The theorem will be proven if we can show that there exists a
$K_{crit}$, independent of $\mathcal{D}_4'$, so that if $|\kb|\geq K_0
> K_{crit}$ then 
\begin{equation}
  \label{3dInward1}
  \Bigl| 2 \pi \sum_{\substack{l_1+l_2=\kb\\ 
      l_1 , l_2 \in \mathcal{Z}}}   \Big[\bigl(
  u^{(1)}_{l_1}(t),\kb\bigr)\omega^{(2)}_{l_2}(t)+
  \bigl(u^{(2)}_{l_1}(t),\kb\bigr)\omega^{(1)}_{l_2}(t) - \bigl(
  \omega^{(2)}_{l_1}(t),\kb\bigr)u^{(1)}_{l_2}(t) -
  \bigl( \omega^{(1)}_{l_1}(t),\kb\bigr)u^{(2)}_{l_2}(t)\Big]\Bigr| 
  <4 \pi^2 \nu |\kb|^\alpha \omega_{\kb}^{(1)}(t)  \ . \notag
\end{equation}
Other boundary points have the same structure so we will only show the
details of the calculation for this case. 

We need to estimate the summation. The total sum is made of smaller
sums which are dominated by sums of the form $\sum_{l_1+l_2=\kb}
|u_{l_1}^{(a)}||\omega_{l_2}^{(b)}||k|$ with $a,b \in \{1,2\}$. As before, we split this sum into three
parts:
\begin{align*}
  \Sigma_1=&\sum_{|l_1|\leq |\frac{\kb}{2}| }
  |u_{l_1}^{(a)}||\omega_{l_2}^{(b)}||k|\\
 \Sigma_2=&\sum_{ |\frac{\kb}{2}| < |l_1|\leq 2|\kb| }
 |u_{l_1}^{(a)}||\omega_{l_2}^{(b)}||k|\\
 \Sigma_3=&\sum_{2|\kb| < |l_1| } |u_{l_1}^{(a)}||\omega_{l_2}^{(b)}||k|
\end{align*}
In $\Sigma_1$, $|l_2| \geq |\frac{\kb}{2}|$ and hence
\begin{align*}
  \Sigma_1 &\leq \frac{\mathcal{D}_4'}{|\kb|^r}2^r |\kb|
  \Biggl(\sum_{|l_1|\leq |\frac{\kb}{2}|} |u_{l_1}^{(a)}|^2
  \Biggr)^{\frac{1}{2}} \Biggl(\sum_{|l_1|\leq |\frac{\kb}{2}|} 1
  \Biggr)^{\frac{1}{2}}\\
& \leq \frac{\mathcal{D}_4'}{|\kb|^r} 2^r \mbox{(const)} \sqrt{E^*} |\kb|^\frac{5}{2}
\end{align*}
The constant is defined by 
\begin{align*}
  \Big(\sum_{|l_1|\leq |\frac{\kb}{2}|} 1\Big) \leq  \mbox{(const)}^2
  |\kb|^3 \ .
\end{align*}
For $\Sigma_2$, we know that $|l_2| \leq 3|\kb|$ and $|u_{l_1}^{(a)}|
\leq \frac{\mathcal{D}_4'}{|l_1|^{r+1}}$
which gives
\begin{align*}
   \Sigma_2 & \leq \frac{\mathcal{D}_4'}{|\kb|^r} 2^r \frac{2}{|k|}
   |k| \sum_{|l_2|\leq 3|\kb|} |\omega_{l_2}^{(b)}| 
  \leq  \frac{\mathcal{D}_4'}{|\kb|^r} 2^{r+1}  \sum_{|l_2|\leq
     3|\kb|} |l_2||u_{l_2}^{(a)}| \\
& \leq  \frac{\mathcal{D}_4'}{|\kb|^r} 2^{r+1}  3|\kb|  \Biggl(\sum_{|l_2|\leq 3|\kb|} |u_{l_1}^{(a)}|^2
  \Biggr)^{\frac{1}{2}} \Biggl(\sum_{|l_2|\leq 3|\kb|} 1
  \Biggr)^{\frac{1}{2}} \\
&\leq    \frac{\mathcal{D}_4'}{|\kb|^r} 2^{r+1}  3   \mbox{(const)}
\sqrt{E^*}|\kb|^{\frac{5}{2}} \ .
\end{align*}
Here the constant is the analogue of the constant in the estimation of
$\Sigma_1$. For $\Sigma_3$, we know that $|l_2| \geq |\kb|$ and thus
\begin{align*}
  \Sigma_3 &\leq |\kb|  \Biggl(\sum_{|l_1|\geq 2|\kb|} |u_{l_1}|^2
  \Biggr)^{\frac{1}{2}} \Biggl(\sum_{|l_2|\geq|\kb|} |\omega_{l_2}|^2 \Biggr)^{\frac{1}{2}}
  \leq |\kb|   \mathcal{D}_4'  \sqrt{E^*} \Biggl(\sum_{|l_2|\geq|\kb|}
  \frac{1}{|l_2|^{2r}} \Biggr)^{\frac{1}{2}}\\
  &\leq  |\kb|   \mathcal{D}_4'  \sqrt{E^*}
  \frac{(\mbox{const})}{|\kb|^{r-\frac{3}{2}}} \leq
  \frac{\mathcal{D}_4'}{|\kb|^r} (\mbox{const})  \sqrt{E^*}
  |\kb|^{\frac{5}{2}} \ .
 \end{align*}
Collecting the three estimates together we see that there is a constant,
depending only on $r$, so that 
\begin{align}
  \label{thm5sumEst}
  \sum_{l_1+l_2=\kb} |u_{l_1}^{(a)}||\omega_{l_2}^{(b)}||k| \leq
  (\mbox{const}) \frac{\mathcal{D}_4'}{|\kb|^r} \sqrt{E^*}
  |\kb|^{\frac{5}{2}} \ .
\end{align}
Using this estimate, we see that the condition in \eqref{3dInward1}
will hold if 
\begin{align*}
  8 \pi  (\mbox{const}) \sqrt{E^*}
  |\kb|^{\frac{5}{2}} \frac{\mathcal{D}_4'}{|\kb|^r}< 4 \pi^2 \nu
  |\kb|^\alpha \frac{\mathcal{D}_8'}{|\kb|^r} \ .
\end{align*}
Since $\alpha > \frac{5}{2}$, this will hold for all $\kb$
sufficiently large; this sets the level of $K_{crit}$. Notice that it
does not depend on $\mathcal{D}_8'$ as was required.
\end{proof}

\begin{proof}[Proof of Theorem \ref{3DifOnlyItWereTrue}]
  The proof of this theorem is similar to the proof of theorem
  \ref{thm:algebraicExponential}. From the assumptions, we know that $
  |\omega_k(t)| \leq \sqrt{\sum_{\ZZ^3} |\omega_l(t)|^2} < \sqrt{\ES^*}$
  for all $t>0$.  We set $a_k^{(j)}(t)=u_k^{(j)} e^{\gamma_7 t |k|}$
  and $b_k^{(j)}(t)=\omega_k^{(j)} e^{\gamma_7 t |k|}$ for $j=1,2$,
  where $\gamma_7$ is a constant we will set later.
  
  Set $\mathcal{D}_7'=2 \max(\sqrt{\ES^*},\mathcal{D}_7)$. Fixing a
  $K_0$, choose $\gamma_7(K_0)$ so that for all $t \in [0,t_1]$, $j \in
  \{1,2\}$, and $k$ with $|k|\leq K_0$, one has $|b^{(j)}_k(t)| \leq
  \frac{\mathcal{D}_7'}{|k|^r}$ .  Notice that by the assumption on
  the initial conditions, we have $|b^{(j)}_k(0)| \leq
  \frac{\mathcal{D}_7'}{|k|^r}$ for all $k$. Consider the set,
  \begin{align*}
    \Omega_7(K_0)=\left\{\left(b_k^{(1)},b_k^{(2)}\right)_{k\in
  \ZZ^2\backslash 0} \mbox{ with } |b_k^{(j)}| \leq
  \frac{\mathcal{D}_7'}{|k|^r} \mbox{ for } j=1,2 \mbox{ and } |k|
  > K_0 \right\} \ . 
  \end{align*}
  As before we will show that, if $K_0$ is chosen large enough, a
  point starting in $\Omega_7$ cannot leave  $\Omega_7$ because the vector
  field along $\partial \Omega_7$ is pointing inward.
  
  We pick a boundary point. For simplicity, we pick the point where
  $b_{\kb}^{(1)}= \frac{\mathcal{D}_7'}{|k|^r}$ and all other
  variables satisfy the inequalities defining $\Omega_7$. 
  In terms of the new variables, the relevant equation of motion reads
  \begin{align*}
    \frac{d b_k^{(1)}(t)}{dt}& = (\gamma_7 |k| -4 \pi^2
    \nu |k|^\alpha) b_k^{(1)}(t) -   2 \pi \sum_{\substack{l_1+l_2=k\\ 
        l_1 , l_2 \in \mathcal{Z}}} \Bigl[ \bigl(
    a^{(1)}_{l_1}(t),k\bigr)b^{(2)}_{l_2}(t) \\ &+
    \bigl(a^{(2)}_{l_1}(t),k\bigr)b^{(1)}_{l_2}(t) - \bigl(
    b^{(2)}_{l_1}(t),k\bigr)a^{(1)}_{l_2}(t) -
    \bigl( b^{(1)}_{l_1}(t),k\bigr)a^{(2)}_{l_2}(t) \Bigr]\frac{e^{-\gamma_7 t
      |l_1|}e^{-\gamma_7 t |l_2|}}{e^{-\gamma_7 t |k|}} \ .
    \end{align*}
  Since $|a^{(j)}_k(t)|=\frac{|b^{(j)}_k(t)|}{|k|}$, to insure that the vector
  field points inward it is sufficient to show that
  \begin{align*}
    2\pi \sum
      |b^{(1)}_{l_1}||b^{(2)}_{l_2}|\frac{|\kb|}{|l_1|}+ 
      |b^{(2)}_{l_1}||b^{(1)}_{l_2}| \frac{|\kb|}{|l_1|} +
      |b^{(2)}_{l_1}||b^{(1)}_{l_2}|\frac{|\kb|}{|l_2|} &
      +|b^{(1)}_{l_1}||b^{(2)}_{l_2}|\frac{|\kb|}{|l_2|} \\ & <
    (4 \pi^2 \nu |\kb|^\alpha - \gamma_7|\kb|)
    \frac{\mathcal{D}_7'}{|\kb|^r} \ .
  \end{align*}
  Each of the terms in the above sum can be estimated with the aid of
  lemma \ref{l:sum}. This transforms the previous condition into
  \begin{align*}
    8 \pi(\mbox{const}) \left( 2^r|\kb| + 2^{r+1} (6|\kb|+1)^{\frac{3}{2}} +
      \frac{1}{2}|\kb|^{2-r}\right) \frac{(\mathcal{D}_7')^2}{|\kb|^r}<
    (4 \pi^2 \nu |\kb|^\alpha - \gamma_7|\kb|)
    \frac{\mathcal{D}_7'}{|\kb|^r} \ .
  \end{align*}
  By picking $K_0$ large enough, we can force this condition to
  hold. The fact that $\gamma_7$ depends on $K_0$ is not a problem
  since it decreases as $K_0$ increases.
  
  This establishes that the vector field points inward along the
  boundary of $\Omega_7$ for all times in the interval $[0,t_1]$.
  Thus at time $t_1$, the trajectory is still in $\Omega_7$. By
  returning to the original variables, we have the desired estimate at
  time $t_1$.
\end{proof} 

\begin{proof}[Proof of Theorem \ref{3DifOnlyItWereTrue2}]
  The proof of this theorem begins as the above theorem and then
  proceeds as the proof of theorem \ref{thm:exponentialExponential}.
  We change variables to $a_k^{(j)}(t)=u_k^{(j)} e^{\gamma_8 |k|}$
  and $b_k^{(j)}(t)=\omega_k^{(j)} e^{\gamma_8 |k|}$. We use the
  assumption on $\sqrt{\sum_{\ZZ^3} |\omega_l(t)|^2}$ to control the
  lower modes. Then we use the estimates from lemma \ref{l:sum} to
  control the nonlinearity. We omit the details.
\end{proof}

\section{Acknowledgements}
The authors thank W.E., C. Fefferman, U. Frisch, G. Gallavotti, J. Mather, V. I. Judovich, F.
Planchon, P. Sarnak, T. Spencer, T. Suidan, J. Vinson, and V. Yakhot for useful
discussions.  The second author thanks NSF (grant DMS-97067994) for
financial support.

\end{document}